\documentclass[11pt,twoside]{article}

\usepackage{epsfig,amsfonts,color}
\usepackage{amsmath}
\usepackage{amssymb, palatino, geometry,url}
\usepackage{algorithmic}
\usepackage[noresetcount,lined,boxed]{algorithm2e} 

\usepackage[colorlinks=true,linkcolor=blue,citecolor=blue,urlcolor=blue]{hyperref}
\usepackage{placeins}

\usepackage{fancyhdr}
\newcommand{\be}{\begin{equation}}
\newcommand{\ee}{\end{equation}}
\newcommand{\bea}{\begin{eqnarray}}
\newcommand{\eea}{\end{eqnarray}}

\newcommand{\bvec}{\left(\begin{array}{c}}
\newcommand{\evec}{\end{array}\right)}
\newcommand{\bsub}{\begin{subequations}}
\newcommand{\esub}{\end{subequations}}

\usepackage{lineno}

\usepackage{tabularx} 
\usepackage{svg} 
\usepackage{makecell} 
\usepackage{amsthm} 
\usepackage{verbatim} 
\usepackage{xcolor}
\usepackage{siunitx}

\theoremstyle{definition}

\theoremstyle{remark}

\theoremstyle{corollary}

\theoremstyle{lemma}

\begin{document}

\title{Recasting AI Data Centers as Engines for Carbon Removal}

\author{
    Zhicong Fang, Boyu Zhang, Jin Shang, Jiaze Ma \\
    \small City University of Hong Kong, Tat Chee Avenue, Kowloon Tong, Hong Kong
}

\date{}

\renewcommand{\thefootnote}{\arabic{footnote}} 

\maketitle

\begin{abstract}

AI data centers (AIDCs) are rapidly increasing electricity demand and associated CO$_2$ emissions, yet they also generate continuous low-grade waste heat. Here we assess whether this heat can be upgraded by heat pumps to drive direct air capture (DAC) and reduce the climate impact of AI infrastructure. We develop a thermodynamically integrated DAC--AIDC system and conduct a region-resolved assessment across the United States, accounting for AIDC capacity, server composition, local climate, electricity prices, and grid carbon intensity. We find that AIDC waste heat can substantially improve net CO$_2$ removal and lower the levelized cost of capture. In carbon-intensive regions, integration can flip DAC from net-positive to net-negative. Under a 2030 scenario with more GPU-intensive AIDCs and cleaner grids, several states achieve removal ratios above 1, indicating that integrated systems can offset their own operational emissions and deliver additional carbon removal. 

\end{abstract}

{\bf Keywords}: AI Data Centers, Direct Air Capture, Thermal Integration\\
\section*{Introduction}
As the physical backbone of artificial intelligence, AI data centers (AIDCs) are expanding at an unprecedented pace, accompanied by a rapid escalation in electricity demand. This rapid growth should also be interpreted in light of historical efficiency gains in IT and cooling infrastructure of DCs, which have strongly shaped the relationship between service demand and electricity use \cite{masanet2020recalibrating}. In 2024, global DC electricity consumption reached approximately 415~TWh and is projected to rise to 945~TWh by 2030—nearly equivalent to the annual electricity consumption of Japan~\cite{IEA2025_energy_ai}. In the United States alone, DCs accounted for 4.4\% of total electricity consumption in 2023, with projections suggesting an increase to 6.7--12\% by 2028~\cite{shehabi20242024}. This surge is driven by the proliferation of compute-intensive AI applications, particularly large language models, and intensified global competition in AI development \cite{de2023growing}. While AI has become a powerful engine for productivity and innovation, it has simultaneously introduced a rapidly growing carbon footprint. Training a single large-scale model such as GPT-3 has been estimated to emit 552~t CO$_2$-equivalent~\cite{patterson2021carbon}, and a bottom-up assessment indicates that 2{,}132 DCs in the United States emitted over 105~Mt CO$_2$-equivalent in 2023, accounting for 2.18\% of national emissions~\cite{guidi2024environmental}. This level of emissions is comparable to the annual tailpipe emissions of more than 22~million passenger vehicles~\cite{epa2023_passenger_vehicle_emissions}. Looking forward, continued expansion of AI server infrastructure is expected to contribute an additional 24--44~Mt CO$_2$-equivalent between 2024 and 2030~\cite{xiao2025environmental}. These trends highlight a fundamental tension: the very infrastructure enabling AI-driven decarbonization across sectors is itself emerging as a significant and rapidly growing source of emissions.

A defining yet underexploited characteristic of DCs is their thermodynamic nature: nearly all electrical energy consumed is ultimately dissipated as heat, which in turn requires additional electricity for removal. Cooling systems already account for 30--40\% of total DC electricity consumption~\cite{zhang2023global}, and this burden is intensifying as rack power densities increase from 8~kW to 17~kW within just two years, with hyperscale GPU clusters exceeding 50~kW. This trend is driving a rapid transition from conventional air cooling to liquid cooling architectures~\cite{azarifar2024liquid,mckinsey2024_ai_power}. The coolant is the mixture of deionized water and anti-freezing solution \cite{alissa2025using}. Consequently, modern AIDCs generate large volumes of low-grade waste heat, typically in the form of 50--60~$^\circ$C return deionized water~\cite{ebrahimi2014review}, representing a substantial yet underutilized thermal resource. Existing efforts primarily focus on integrating this waste heat into district heating systems, where notable success has been demonstrated—for example, achieving up to 36\% CO$_2$ reduction in Finland and supplying up to 80\% of local heating demand in specific deployments~\cite{tervo2025reducing,townsend2024our}. However, such approaches are inherently constrained by geographic proximity, climate conditions, and the limited transportability of low-grade heat~\cite{santin2020feasibility,jing2026large}. Moreover, heating demand is highly seasonal and mismatched with the continuous operation of AIDCs. These limitations underscore the need for more generalizable and scalable pathways to valorize AIDC waste heat. Developing such approaches could simultaneously reduce the parasitic energy demand of cooling systems and unlock a large, low-cost thermal source for integrated decarbonization technologies.

\begin{figure}[htbp]
    \centering
    \includegraphics[width=1\linewidth]{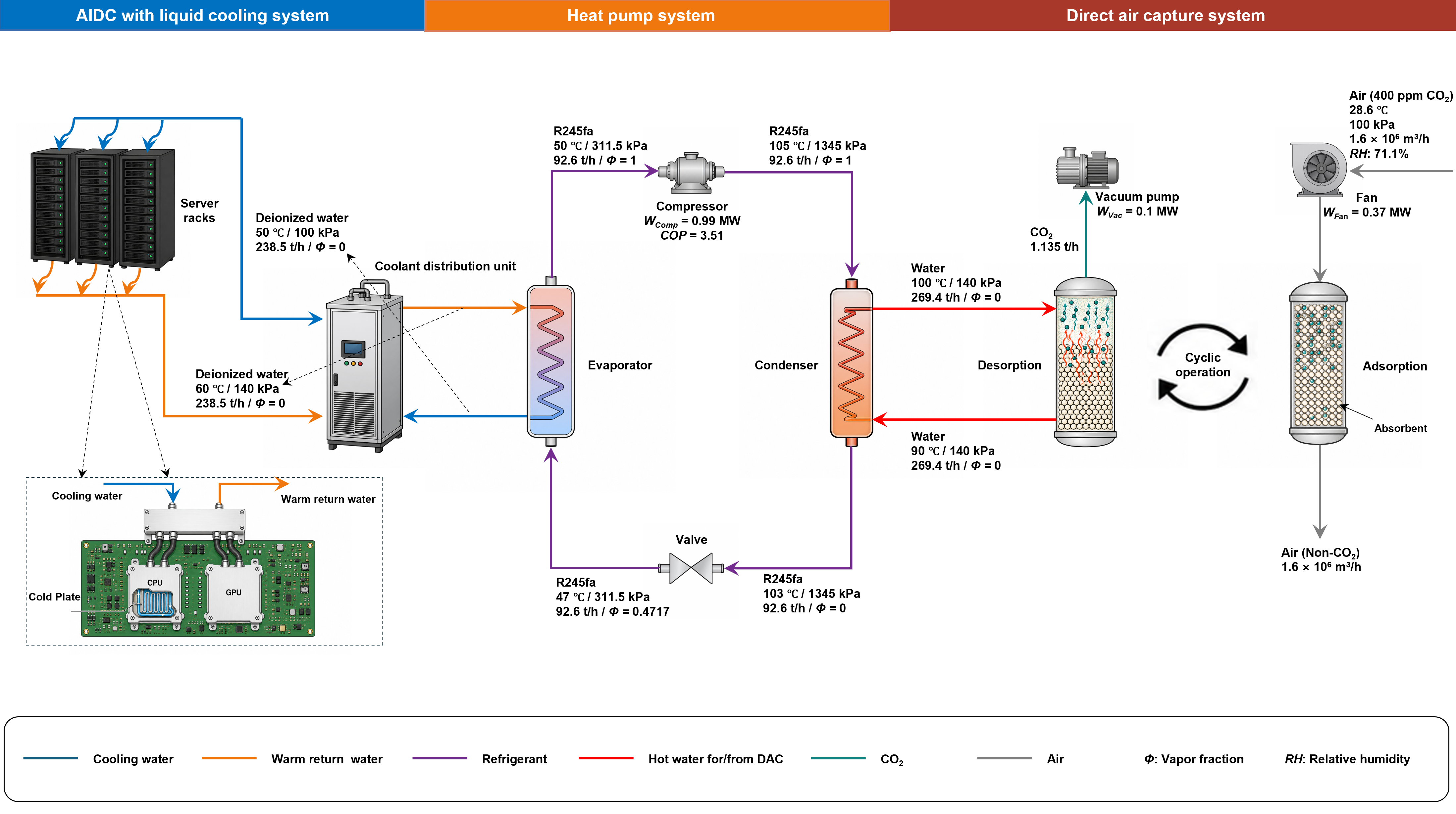}
    \caption{Thermally integrated DAC--AIDC system concept. Low-grade waste heat from an AIDC liquid-cooling loop (50~$^\circ$C supply, 60~$^\circ$C return) is upgraded by a heat pump to provide $\sim$100~$^\circ$C regeneration heat for a temperature--vacuum swing adsorption (TVSA) DAC unit. Ambient air is driven through absorbent beds for CO$_2$ capture during adsorption, followed by vacuum/thermal regeneration to release a concentrated CO$_2$ stream for utilization or sequestration; auxiliary electricity is required for the compressor, fans, and vacuum pump, linking net removal performance to local grid conditions.}
    \label{fig1 DACDC_process}
\end{figure}

Among negative-emissions technologies, direct air capture (DAC) stands out for its scalability and location flexibility. As illustrated in Fig.~\ref{fig1 DACDC_process}, a typical DAC system operates through cyclic adsorption and desorption processes: ambient air is drawn through contactors where CO$_2$ is selectively captured by solid absorbents, followed by thermal regeneration to release a concentrated CO$_2$ stream for utilization or storage. Unlike post-combustion capture, DAC is decoupled from point emission sources and land-use constraints~\cite{keith2018process,mcqueen2021review,fout2022direct}, making it a great option for achieving deep decarbonization. Under the IEA Net-Zero pathway, DAC is expected to remove approximately 980~Mt CO$_2$ annually by 2050~\cite{international2022direct}. 
Despite these advantages, large-scale deployment of DAC is fundamentally constrained by its substantial thermal energy demand. In adsorption-based DAC systems, solid absorbents capture CO$_2$ from ambient air and require regeneration at elevated temperatures, typically in the range of 80--120~$^\circ$C~\cite{zhou2021low_temp_dac}. This regeneration step constitutes the primary energy bottleneck. State-of-the-art temperature--vacuum swing adsorption (TVSA) systems require 4--6~GJ of thermal energy and 0.5--1.5~GJ of per tonne of CO$_2$ captured~\cite{sabatino2021comparative,sievert2024considering,wiegner2022optimal}. The high capital investment of DAC system and the high operational cost due to massive energy consumption 
cause unaffordable levelized cost of carbon capture (LCCC). Standalone DAC with electric heating exhibits LCCC around 773--802~\$/t \cite{patel2024techno}. Energy price is critical factor for operational cost. When electricity price rises from 50~\$/MWh to 100~\$/MWh, case study of standalone DAC in US reveals the LCCC increasing from 320--540~\$/t to 405--659~\$/t \cite{sendi2022geospatial}. To alleviate this challenge, waste heat recovery has been explored as a promising pathway \cite{d2024integrating}. For example, coupling DAC with industrial waste heat in refinery systems has reduced the LCCC to 148.5~\$/t after large-scale deployment in the future ~\cite{odeh2025techno}. However, such heat sources are geographically constrained, progressively diminished under industrial decarbonization, and often insufficient in temperature. To bridge the temperature gap between low grade industrial waste heat and the DAC, heat pumps offer a viable solution. When integrated, heat pumps can reduce DAC energy consumption by up to 69.5\%~\cite{ge2024innovative}, and can even enable the utilization of ambient heat under advanced configurations~\cite{leonzio2022innovative}. 

In this context, the low-grade waste heat (e.g., $\sim$50--60$^\circ$C) from AI data centers (AIDCs) presents a particularly attractive and rapidly growing thermal resource that can be upgraded via heat pumps to supply DAC systems~\cite{diaz2025flipping}. However, the overall energy efficiency and carbon abatement performance of DAC remain highly sensitive to local conditions, including ambient temperature, humidity, grid carbon intensity, and the availability of waste heat from AIDCs. Maximizing the carbon capture benefit therefore requires careful co-design and optimization of system configuration and operation under region-specific conditions. The key unanswered question, and also  the central focus of this work, is \textit{whether the waste heat from AIDCs can be effectively harnessed via heat pump integration to drive DAC systems, thereby reconciling the rapid expansion of AI infrastructure with climate mitigation goals}. Given the substantial spatial heterogeneity in climate conditions, grid carbon intensity, and waste heat availability, an equally critical question is: what is the maximum achievable carbon abatement potential of DAC systems driven by AIDC waste heat across different regions in the United States?

To quantify the decarbonization potential of harnessing AIDC waste heat for carbon capture, we propose a thermodynamically integrated DAC--AIDC system in which a heat pump upgrades low-grade heat from AI server cooling loops to the $\sim$100~$^\circ$C required for absorbent regeneration (Fig.~\ref{fig1 DACDC_process}). In modern AIDCs, the high heat flux generated by CPUs and GPUs is removed via liquid cooling loops operating at approximately 50~$^\circ$C supply and 60~$^\circ$C return~\cite{ebrahimi2014review}, providing a continuous and concentrated source of low-grade thermal energy. This heat is upgraded through a heat pump cycle and supplied to a temperature--vacuum swing adsorption (TVSA)-based DAC process, where CO$_2$ is captured from ambient air and subsequently released during thermal regeneration for utilization or sequestration.

Despite this thermodynamic synergy, DAC operation remains inherently electricity-dependent because compressors, fans, and auxiliary units must be powered by the grid. Consequently, the net carbon abatement and economic performance of the integrated system are strongly coupled to local grid carbon intensity and electricity prices. Moreover, ambient temperature and humidity directly influence adsorption capacity and regeneration heat demand, introducing pronounced geographic variability in capture efficiency. To ensure a realistic assessment, our nationwide analysis explicitly accounts for (i) the spatial distribution and capacity of AIDCs, (ii) the heterogeneous mix of AI server types (e.g., GPU-accelerated versus non-accelerated configurations with different rated powers and waste-heat profiles), and (iii) local climate and grid conditions. These coupled dependencies make DAC--AIDC viability fundamentally location-specific, motivating the region-resolved evaluation of waste heat availability, net carbon removal, and capture cost presented in the following sections.

Based on a nationwide assessment across the United States, we show that fully utilizing AIDC waste heat in 2026 could enable up to 25.6~MtCO$_2$/yr of carbon capture, comparable to the annual emissions of a small developed country. More importantly, the integration (Fig.~\ref{fig1 DACDC_process}) unlocks additional carbon removal potential by shifting regions with high grid carbon intensity from net-positive to net-negative outcomes--- effectively \emph{flipping} the carbon abatement feasibility in locations where standalone DAC fails to deliver net removal. In several regions, we further find that AIDC waste heat can drive DAC operation to fully offset the indirect CO$_2$ emissions associated with AIDC electricity use, turning AIDCs from emission liabilities into net removal assets.

\section*{Results}

\subsection*{Conceptual process description and thermodynamic analysis}
Liquid cooling in AIDCs basically includes cold-plate cooling and immersion cooling \cite{alissa2025using}. The former accounts for the vast majority of liquid-cooled capacity and is projected to take the lead through the end of the decade \cite{DellOroGroup2026LiquidCoolingMarket,Davis2024WaterColdPlatesDLC}. Since we focus on the cold-plate cooling in this paper. The core integration challenge is a temperature--grade mismatch: AIDCs reject abundant but \emph{low-grade} heat through liquid cooling (typically 50~$^\circ$C supply and 60~$^\circ$C return), whereas adsorption-based DAC regeneration requires a \emph{higher-grade} thermal stream near 100~$^\circ$C. We design a single-stage vapor-compression heat pump cycle that extracts heat from the AIDC return-water loop in the evaporator and delivers heat to the DAC hot-water loop in the condenser. The cycle is configured to match the AIDC waste-heat temperature window and minimize electricity consumption while respecting realistic heat-exchanger temperature approaches. R245fa is selected as the working fluid because its saturation properties and pressure levels are well suited for upgrading heat from the 50--60~$^\circ$C range to $\sim$100~$^\circ$C in a compact, single-stage configuration (details in SI). A high-fidelity process simulation is performed in Aspen Plus, with all parameter values and mass/energy balances reported in the Supplementary Information.

\begin{figure}
    \centering
    \includegraphics[width=0.75\linewidth]{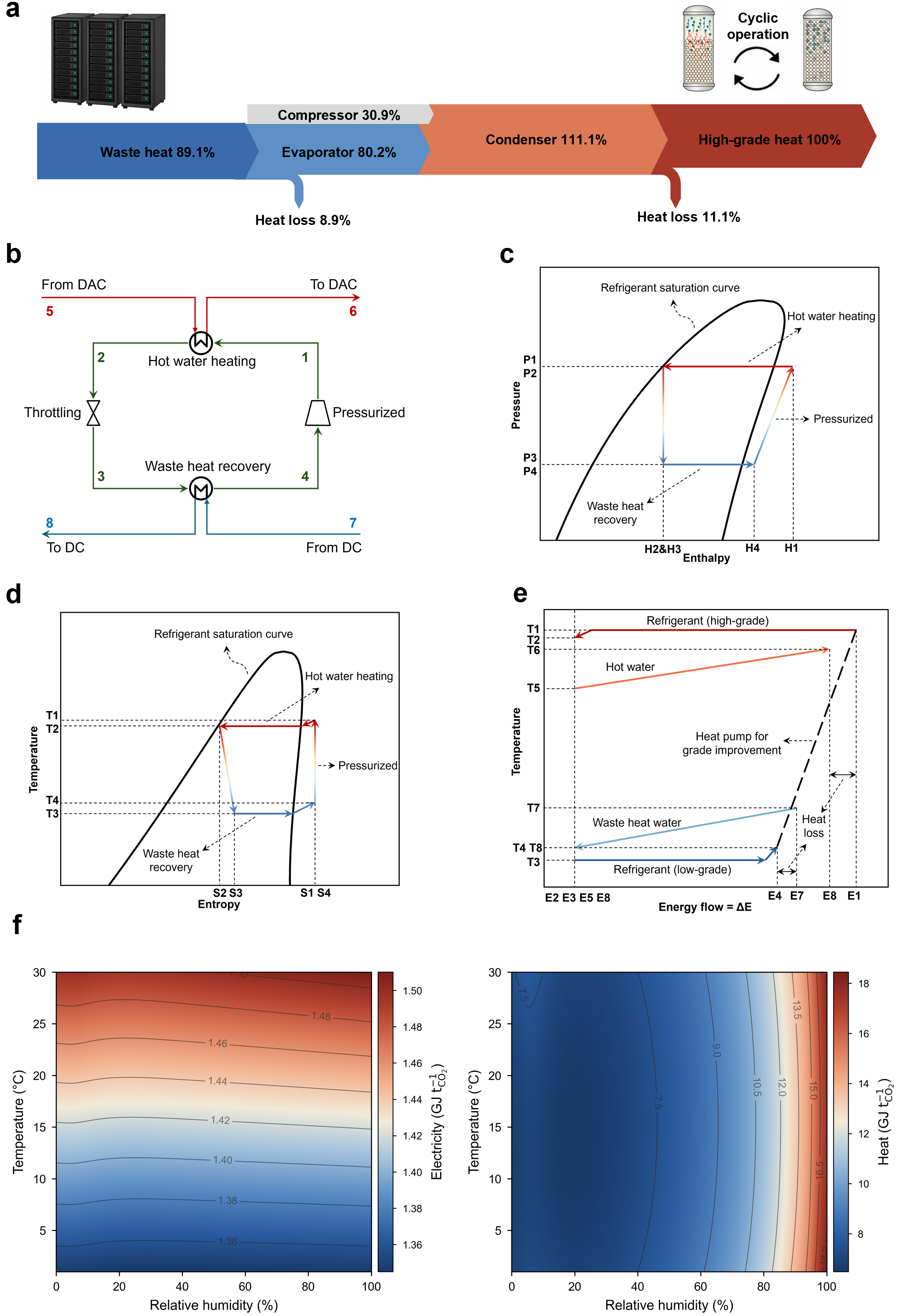}
    \caption{Thermodynamic rationale and performance of the heat-pump-enabled DAC--AIDC integration. (a) System energy flows and losses (Sankey). (b) Vapor-compression heat pump configuration for upgrading AIDC waste heat to DAC regeneration heat. (c) Pressure--enthalpy (P--$h$) diagram and (d) temperature--entropy (T--$s$) diagram illustrating the upgrade mechanism. (e) Temperature-level (composite-curve) check showing required driving forces across evaporator and condenser. (f) Electricity and regeneration heat demand per tonne CO$_2$ captured under varying ambient temperature and relative humidity.}
    \label{thermodynamics}
\end{figure}

The heat pump operates as a standard vapor-compression cycle (Fig.~\ref{thermodynamics}b), which can be read in four steps: 
(i) \emph{Evaporation:} low-pressure refrigerant absorbs heat from the AIDC return water and vaporizes; 
(ii) \emph{Compression:} electricity drives a compressor that raises the vapor pressure and temperature; 
(iii) \emph{Condensation:} the hot, high-pressure refrigerant releases heat to the DAC hot-water loop, delivering $\sim$100~$^\circ$C regeneration heat; and 
(iv) \emph{Expansion:} throttling reduces the pressure to close the cycle. 
The P--$h$ and T--$s$ diagrams (Fig.~\ref{thermodynamics}c--d) visualize the same mechanism: compressor work increases the refrigerant temperature so that heat can be rejected at a higher temperature than the original 50--60~$^\circ$C waste-heat source.

Fig.~\ref{thermodynamics}a summarizes the energy balance. Accounting for heat-exchanger losses of $\sim$10\%, the upgraded waste heat can fully meet the DAC regeneration duty. Supplying 1~kJ of useful regeneration heat requires $\sim$0.30~kJ of electricity, corresponding to $\mathrm{COP}=3.51$. This relatively high COP is enabled by the elevated-temperature heat source (AIDC return water). If ambient air is used instead, the COP drops to $\sim$2~\cite{leonzio2022innovative}, implying a substantially larger electricity penalty for the same regeneration heat duty. Fig.~\ref{thermodynamics}e provides a temperature-level check that explains why the integration is feasible. The design maintains (i) refrigerant evaporation temperature below the AIDC return-water temperature to ensure heat uptake and (ii) refrigerant condensation temperature above the DAC hot-water temperature to ensure heat delivery, thus preserving positive driving forces in both heat exchangers. In hot and humid regions where AIDC cooling may require chillers, the same heat pump can partially substitute chiller duty while upgrading heat for DAC regeneration, creating an additional electricity-saving pathway (quantified in the SI). Finally, DAC energy demand remains climate-sensitive: higher humidity increases co-adsorbed water and regeneration heat demand, while temperature affects adsorption behavior and fan electricity through air-density changes. As shown in Fig.~\ref{thermodynamics}f, these effects motivate the spatially resolved assessment in the following sections, where net removal and economics depend on local weather and grid conditions.

\subsection*{Matching AI waste heat with DAC energy demand}

Before assessing the carbon abatement potential of the proposed system, we first quantify the waste heat availability from AIDCs across the United States. We combine the 2024 Lawrence Berkeley National Laboratory report on the installed base and rated power of four AI server classes (AI-8GPU, AI-4GPU, AI-2GPU, and AI Non-accelerated)~\cite{shehabi20242024} with county-level DC capacity from the DOE ``Speed to Power'' initiative~\cite{DOE_SpeedToPower_2025} and state-level electricity consumption estimates from EPRI~\cite{intelligence2024analyzing}. servers have different waste heat structure. The recycled waste heat from AIDC refers to the heat released from processors ( CPUs/GPUs), and memory chips as their upper limit of working temperature is higher than 60~$^\circ$C and require liquid cooling. Other equipment operating at low temperature like AC/DC conversion, disk drives and motherboards is cooled by air. In this study, we assume that this heat accounts for 46\% of total waste heat in AI Non-accelerated servers while this value increases to 89\% for AI-8GPU servers \cite{ahmed2021review,nvidia2025dgxsuperpod}.  The corresponding ratios for AI-2GPU and AI-4GPU servers are then obtained by interpolation based on the number of GPUs, which are 69\% and 78\%, respectively. AI servers are allocated in proportion to each state’s share of the top thirty DC loads and then distributed within each state according to county-level DC capacity. Full assumptions and calculations are provided in the SI.

\begin{figure}
    \centering
    \includegraphics[width=0.65\linewidth]{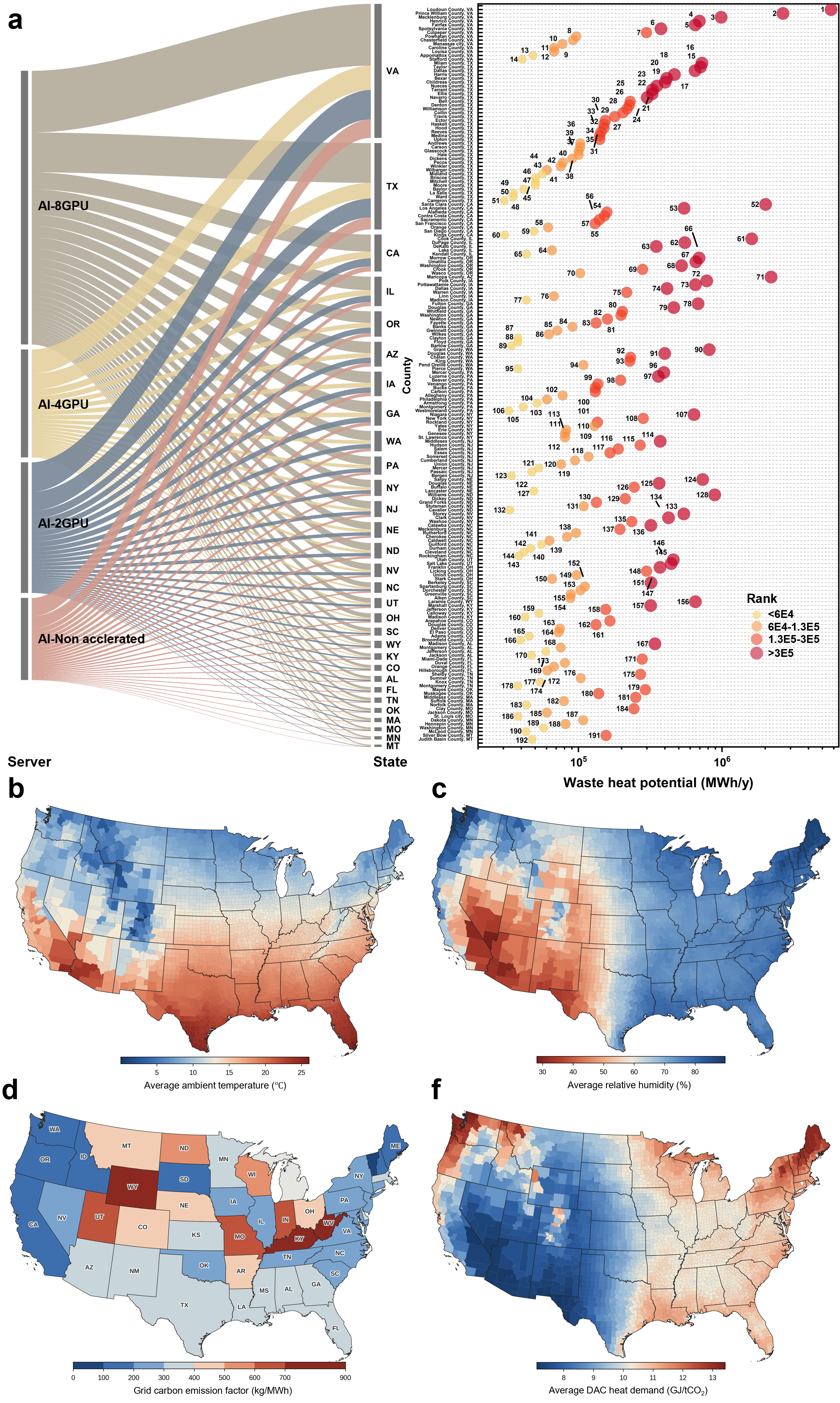}
    \caption{ Spatial heterogeneity in AIDC waste heat availability and DAC energy drivers across the United States. 
(a) Estimated annual waste heat potential of AI servers in U.S. data centers in 2024, disaggregated by server class (AI-8GPU, AI-4GPU, AI-2GPU, and non-accelerated) and aggregated by state; the corresponding county-level distribution (right) highlights strong spatial concentration within a small number of high-load counties. 
(b) County-level annual mean ambient temperature and (c) annual mean relative humidity used to parameterize climate-dependent DAC performance. 
(d) State-level grid carbon emission factor (CEF), which determines indirect CO$_2$ emissions from electricity consumption. 
(f) County-level annual average DAC regeneration heat demand (GJ/t$_{\mathrm{CO_2}}$), showing systematically lower heat demand in hot--dry regions and higher demand in cold--humid regions.
    }
    \label{fig3Wasteheatpotential}
\end{figure}

The resulting waste heat potential exhibits strong spatial heterogeneity and concentration (Fig.~\ref{fig3Wasteheatpotential}a). Virginia and Texas reach 11.9~TWh/y and 7.7~TWh/y, respectively, which are 9.9$\times$ and 6.4$\times$ larger than the 15$^{\text{th}}$ ranked Nevada (1.2~TWh/y). The top 15 states together account for $\sim$85\% of the national total, and the top 8 states each exceed 2~TWh/y. Importantly, the county-level distribution reveals pronounced clustering: although California, Illinois, and Arizona have lower state totals than Texas, their largest counties exhibit substantially higher waste heat potential, reaching $2.0\times10^{6}$~MWh/y (CA), $1.6\times10^{6}$~MWh/y (IL), and $2.2\times10^{6}$~MWh/y (AZ), compared with $7.2\times10^{5}$~MWh/y in the largest Texas county. This hub-like concentration is deployment-relevant: it supports co-located scale-up of DAC and heat-pump infrastructure near a small number of high-load sites, rather than requiring many dispersed small plants.

To translate waste heat availability into carbon removal performance, we also map the key regional drivers of DAC energy intensity and indirect emissions (Fig.~\ref{fig3Wasteheatpotential}b--f). Ambient temperature increases from north to south (Fig.~\ref{fig3Wasteheatpotential}b), while relative humidity is highest in the northwestern and northeastern coastal regions and lowest in the arid southwest (Fig.~\ref{fig3Wasteheatpotential}c). Grid carbon emission factors (CEF) vary substantially across states (Fig.~\ref{fig3Wasteheatpotential}d), directly determining the indirect CO$_2$ emissions associated with electricity consumption. The county-level annual average regeneration heat demand (Fig.~\ref{fig3Wasteheatpotential}f) shows large spatial contrasts and tracks humidity more strongly than temperature: cold--humid regions such as Washington and Oregon exhibit high heat demand (9.2--13.4~GJ/t and 8.6--13.1~GJ/t, respectively), whereas hot--dry states including Arizona, California, and Texas show markedly lower median heat demands (7.3, 9.2, and 9.1~GJ/t). Virginia, despite having the largest waste heat potential, maintains moderate heat demand (median 10.6~GJ/t) due to its comparatively lower humidity and cleaner grid than Texas. 

Taken together, Fig.~\ref{fig3Wasteheatpotential} highlights three coupled determinants of regional DAC--AIDC performance: waste heat availability (capacity), climate-dependent energy intensity (GJ/t), and grid carbon intensity (net removal feasibility). Their spatial patterns are not aligned, motivating the coupled, region-resolved evaluation of net carbon abatement and economics in the following sections.

\subsection*{Net removal ratio and CO$_2$ uptake enabled by DAC--AIDC thermal integration}

\begin{figure}
    \centering
    \includegraphics[width=1\linewidth]{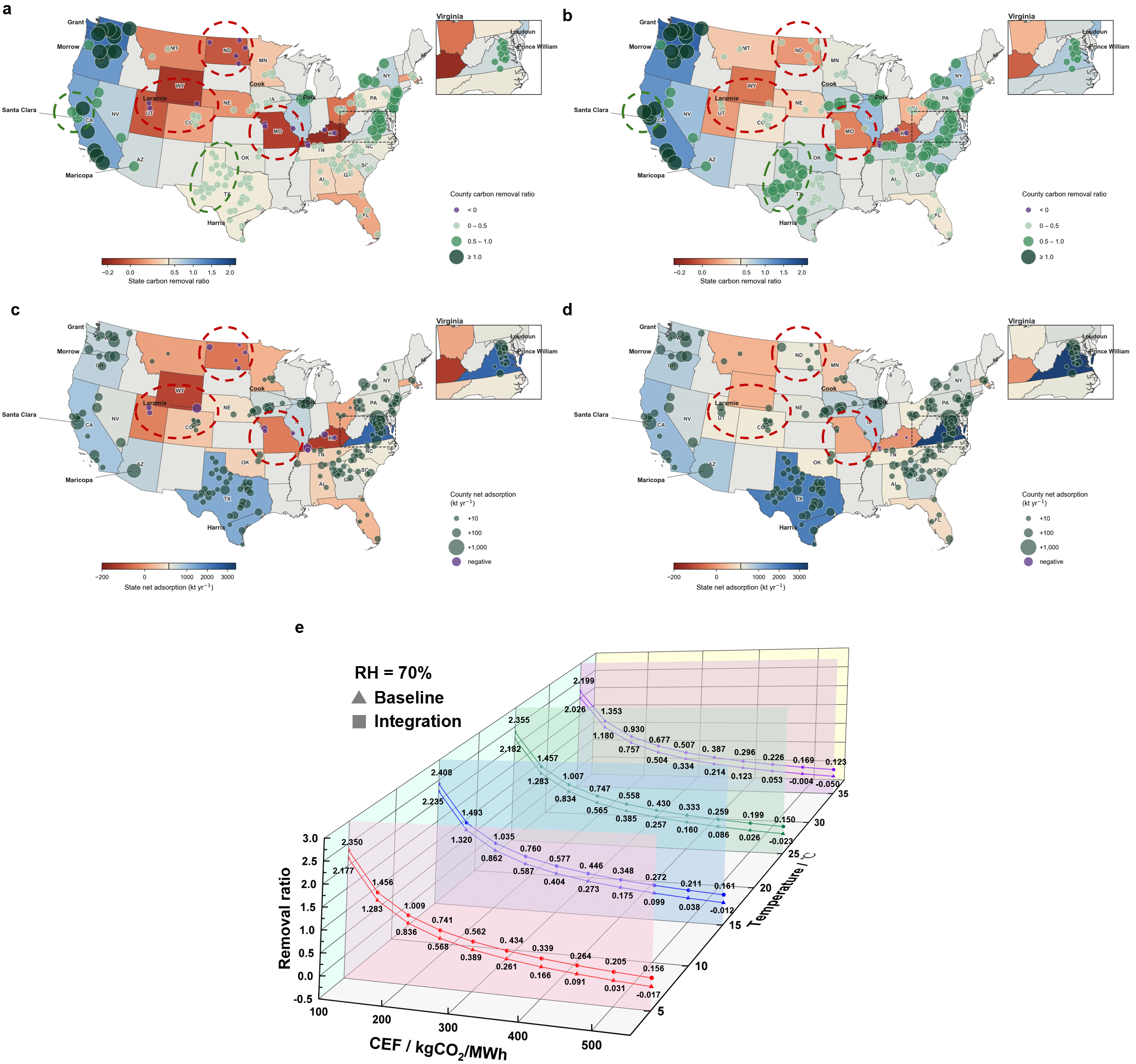}
    \caption{DAC--AIDC integration increases net carbon removal and expands the net-negative deployment space across the United States. (a) Baseline (standalone) DAC net removal ratio, defined as net CO$_2$ removed by DAC after accounting for electricity-related emissions, normalized by AIDC indirect CO$_2$ emissions from electricity use; negative values indicate net-positive emissions and values $>1$ indicate that DAC fully offsets AIDC indirect emissions and removes additional CO$_2$. (b) Net removal ratio under DAC--AIDC integration, showing substantial uplift and ``feasibility flips'' from negative to positive outcomes in multiple high-CEF regions due to waste-heat-driven regeneration. (c,d) Annual CO$_2$ adsorption (uptake) under baseline and integration cases, respectively, quantifying additional capture enabled by AIDC waste heat. (e) Sensitivity of net removal ratio to ambient temperature and grid carbon emission factor (CEF), illustrating that integration shifts the net-zero boundary and improves robustness under carbon-intensive grids.}
    \label{fig4performanceanalysis}
\end{figure}

Fig.~\ref{fig4performanceanalysis} summarizes the central message of this work: \emph{harnessing AIDC waste heat can substantially increase net carbon removal and, critically, flip regions where standalone DAC fails into regions where DAC becomes carbon-negative.} To quantify performance in a way that directly reflects AIDC decarbonization, we define the \emph{net removal ratio} as the net CO$_2$ removed by DAC (gross capture minus electricity-related emissions) divided by the indirect CO$_2$ emissions associated with the AIDC’s electricity consumption. A value above unity is particularly meaningful: it indicates that the integrated system can fully offset the AIDC’s indirect emissions and remove additional CO$_2$, effectively turning the AIDC into a net carbon-negative facility.

Under the baseline scenario (standalone DAC without AIDC waste heat; Fig.~\ref{fig4performanceanalysis}a), more than half of the states exhibit net removal ratios below 0.4, indicating limited capability to offset AIDC-related emissions. Even among the three states with the largest AI server inventories, Virginia and Texas achieve only 0.53 and 0.34 on average, respectively, while California reaches 1.03 due to its cleaner grid. Importantly, baseline deployment can be \emph{counterproductive} in parts of the U.S.: 15 counties across Utah, North Dakota, Missouri, Wyoming, and Kentucky show negative net removal ratios, meaning the DAC system emits more CO$_2$ (via electricity use) than it removes. This failure mode occurs where climate-driven energy demand is high (cold and humid conditions) and grid carbon intensity is unfavorable.

When waste heat is supplied by AIDCs through heat-pump integration (Fig.~\ref{fig4performanceanalysis}b), the performance landscape changes qualitatively. The average net removal ratio increases to 0.71 in Virginia, 0.51 in Texas, and 1.23 in California, corresponding to relative gains of 34\%, 50\%, and 19.4\%, respectively. The most striking outcome is the \emph{feasibility flip}: all previously negative states except Kentucky become net-negative, with county-level ratios shifting from below zero to positive values (up to $\sim$0.2). In other words, waste-heat-driven regeneration converts regions where standalone DAC cannot deliver net removal into regions where DAC becomes a viable abatement option. Meanwhile, multiple counties in clean-grid states achieve ratios above 1, implying complete offsetting of AIDC indirect emissions and additional CO$_2$ removal. This ``$>1$'' regime is a strong indicator of \emph{carbon-negative AI infrastructure} enabled by thermal integration.

Panels (c) and (d) quantify how this flip translates into absolute CO$_2$ removal. Nationally, the annual CO$_2$ adsorption increases from 9{,}751.7~kt/y (baseline) to 13{,}760~kt/y (integration), i.e., an additional 41.1\% capture enabled by AIDC waste heat. Virginia, Texas, and California remain the dominant contributors due to their large AIDC clusters, increasing from 2{,}427.7 to 3{,}202.6~kt/y, 1{,}297.5 to 2{,}023.9~kt/y, and 882.8 to 1{,}044.3~kt/y, respectively. Notably, the integration also unlocks substantial capture in regions that were previously unattractive for deployment: excluding Kentucky (extremely high CEF), the other previously net-positive states collectively achieve 489.2~kt/y of adsorption under integration, demonstrating that waste heat expands the set of regions where meaningful carbon removal is achievable rather than merely improving already-favorable locations.

Fig.~\ref{fig4performanceanalysis}e explains the mechanism behind the flip and identifies when it occurs. Net removal ratio decreases with increasing grid carbon emission factor (CEF) in both scenarios, but integration shifts the curve upward across all operating conditions. At $15~^\circ$C, the baseline net removal ratio drops from 2.24 to $-0.01$ as CEF increases from 100 to 550~kg\,CO$_2$/MWh, whereas the integration case remains positive (2.41 to 0.16) over the same range. Temperature plays a secondary but non-negligible role, becoming decisive near the net-zero boundary: at high CEF (e.g., 500~kg\,CO$_2$/MWh), increasing ambient temperature from 15 to 35~$^\circ$C can shift the baseline case from slightly positive to slightly negative. Crucially, integration provides an approximately constant uplift (about 0.17 in net removal ratio), which is sufficient to push high-CEF regions across the zero threshold. As a result, the feasible deployment window expands substantially: areas that are net-positive under standalone DAC become net-negative under DAC--AIDC integration, consistent with the county-level patterns in Fig.~\ref{fig4performanceanalysis}~(a)--(d).

\subsection*{Economic viability and future removal potential}

\begin{figure}
    \centering
    \includegraphics[width=0.75\linewidth]{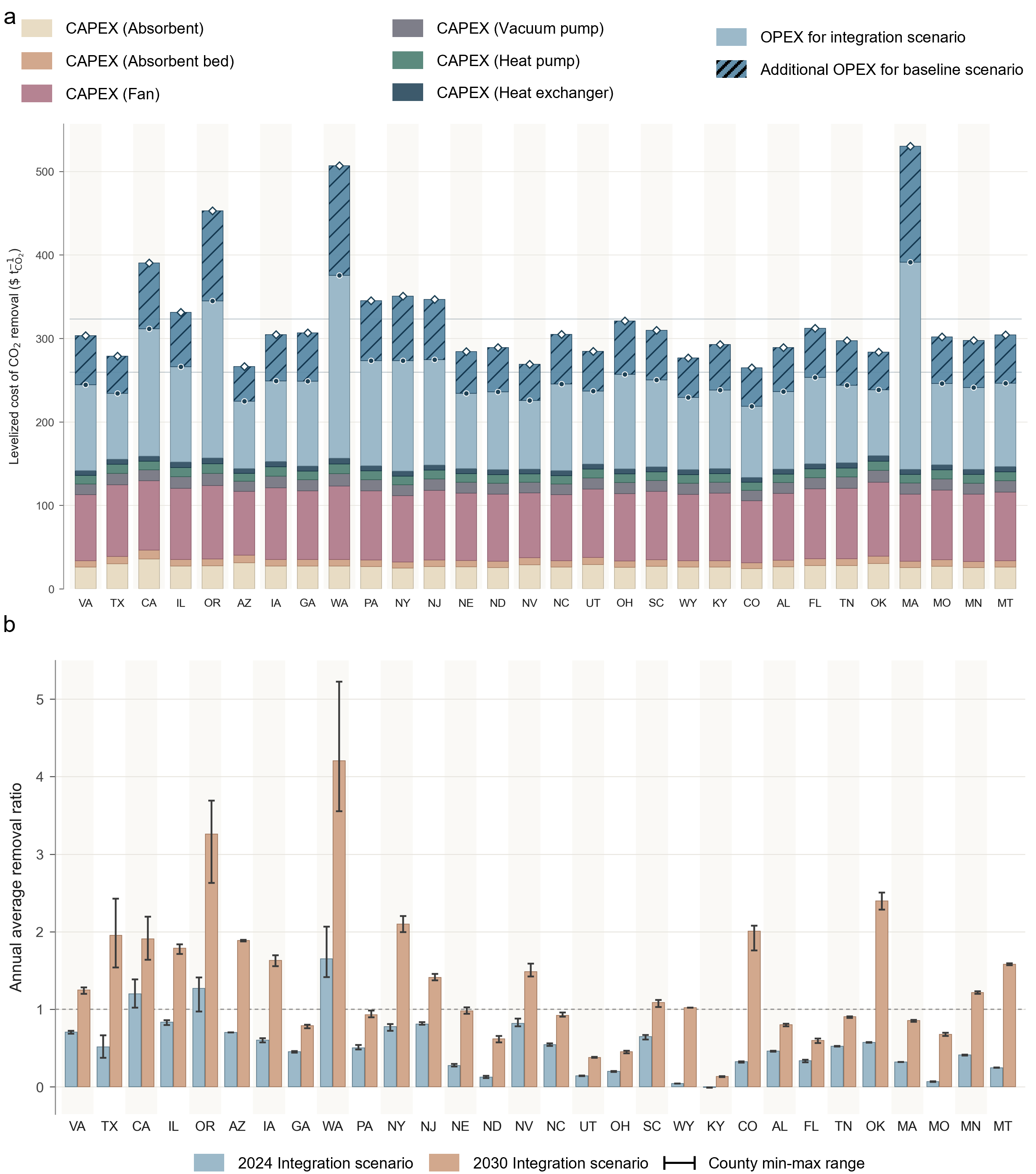}
    \caption{
Economic performance and future carbon-removal potential of DAC--AIDC integration across U.S. states.
(a) Levelized cost of CO$_2$ capture (LCCC) under the integrated DAC--AIDC scenario. The hatched segment indicates the additional OPEX that would be required under the standalone DAC baseline. 
(b) Annual average removal ratio of the integrated DAC--AIDC system in 2024 and 2030. The 2030 scenario assumes a higher share of GPU-based AI computing, which increases recoverable waste heat, together with a cleaner power grid, which reduces indirect emissions from DAC electricity consumption. Error bars indicate county-level minimum--maximum ranges within each state.
}
    \label{fig6LCCC}
\end{figure}

Fig.~\ref{fig6LCCC} shows that DAC--AIDC integration can improve carbon-removal performance while keeping the levelized cost of CO$_2$ capture (LCCC) within a competitive range. The LCCC is mainly driven by electricity-related OPEX; therefore, using AIDC waste heat to reduce DAC energy demand is especially valuable in high-electricity-price states. Under the standalone DAC case, Massachusetts, Washington, Oregon, and California show high LCCC values of 530.4, 506.8, 452.7, and 350.6~\$/t$_{\mathrm{CO_2}}$, respectively. After integration, these costs decrease by approximately 20.0\%--26.2\%, demonstrating that waste-heat utilization can substantially reduce the economic penalty of DAC operation.

The integrated system is also attractive in low-cost electricity markets. Arizona, Texas, and Oklahoma achieve some of the lowest integrated LCCC values, at 224.9, 234.2, and 238.4~\$/t$_{\mathrm{CO_2}}$, respectively. Among CAPEX components, fan cost is the largest contributor after OPEX because DAC requires large air throughput to capture dilute atmospheric CO$_2$. Absorbent replacement is another important cost component, especially in states with stronger seasonal variation in temperature and humidity, where the system must be sized for the least favorable month. These results indicate that favorable deployment regions are not determined by electricity price alone, but by the combined effects of electricity cost, grid-carbon intensity, climate conditions, AIDC capacity, and recoverable waste heat.

Several states already show a strong balance between economic cost and removal performance in 2024. Virginia is particularly important because it hosts large AIDC capacity while maintaining a competitive integrated LCCC, which decreases from 303.4 to 244.5~\$/t$_{\mathrm{CO_2}}$. Texas is also notable: it combines one of the lowest integrated LCCC values with a 2024 removal ratio of approximately 0.5. These cases suggest that DAC--AIDC integration can already provide meaningful carbon-abatement value under current grid and data-center conditions.

The 2030 integrated scenario reveals a much stronger opportunity. According to the estimation by Lawrence Berkeley National Laboratory ~\cite{shehabi20242024}, future AIDCs contain a higher share of 8-GPU servers, of which number will be close to the sum of other types of servers in 2028. This significantly increases the recovered waste heat ratio in AIDCs. The power grid also becomes significantly cleaner predicted in ''Cambium 2024 scenario'' by National Renewable Energy Laboratory \cite{gagnon2025cambium}, which lowers the indirect emissions associated with DAC electricity consumption. These two trends reinforce each other: future AI data centers provide more recoverable heat, while the DAC system becomes less carbon-intensive to operate. As a result, the integrated removal ratio increases sharply from 2024 to 2030 in many states. More details are shown in SI.

The most important finding is that several states achieve integrated removal ratios far above 1 by 2030. This threshold is central to the interpretation of the results. A removal ratio above 1 means that the DAC--AIDC system removes more CO$_2$ than the operational emissions of the associated AIDC. Washington and Oregon reach removal ratios above 4 and 3, respectively, suggesting that one integrated DAC--AIDC system could offset not only its own AIDC emissions, but also the operational emissions of several additional AIDCs of comparable scale. In these regions, DAC--AIDC integration is no longer only an efficiency measure for reducing the footprint of a single facility; it becomes a regional carbon-removal strategy.

This result also has direct implications for future AI infrastructure planning. Regions with favorable climate conditions, cleaner electricity, and high recoverable waste-heat potential should be prioritized for additional AIDC deployment. Washington and Oregon are strong examples because both the power-system carbon intensity and local climate conditions are favorable for net CO$_2$ removal. Texas is another important case: by 2030, its integrated removal ratio approaches 2 while maintaining a low LCCC. This suggests that Texas could become one of the most attractive regions for future DAC-coupled AI infrastructure, as low electricity cost, rapid renewable deployment, and increasing GPU-driven waste heat jointly improve both economic and climate performance. Virginia also exceeds a removal ratio of 1, indicating that integrated DAC can fully compensate for the operational emissions of its coupled AIDCs.

The 2030 results therefore suggest that the siting of future AIDCs should not be based only on land availability, electricity price, and grid interconnection capacity. Carbon-removal potential should also become a siting criterion. Deploying more AIDCs in regions where waste heat can be effectively converted into DAC regeneration energy, and where the power grid is already relatively clean or rapidly decarbonizing, can reduce the carbon burden of AI growth. In the strongest cases, such deployment may even allow AI infrastructure to become a net contributor to regional carbon removal rather than merely a source of additional electricity demand.

Overall, Fig.~\ref{fig6LCCC} shows that DAC--AIDC integration becomes more valuable over time. In 2024, the integrated system already lowers capture cost and improves removal performance in selected states. By 2030, with more GPU-intensive AIDCs and cleaner grids, the integrated system can exceed full AIDC emission compensation in several regions, and in some states remove several times more CO$_2$ than the associated AI infrastructure emits. This provides a clear system-level message: the future decarbonization of AI should combine cleaner electricity, strategic siting, and direct use of data-center waste heat for carbon removal.

\FloatBarrier
\section*{Discussion}

The rapid expansion of AI infrastructure creates a dual challenge for sustainability. On one hand, AI is increasingly positioned as an enabling technology for climate mitigation, scientific discovery, and energy-system optimization. On the other hand, the physical infrastructure required to support AI is becoming a material source of electricity demand and indirect CO$_2$ emissions. This work shows that AI data centers should not be viewed only as electricity-consuming digital infrastructure, but also as spatially concentrated thermal resources. By coupling AIDC waste heat with heat-pump-assisted DAC, the same infrastructure that drives additional power demand can be converted into a platform for engineered carbon removal.

Our results demonstrate that this integration changes the carbon-removal feasibility of DAC in a fundamental way. In regions with carbon-intensive grids, this electricity penalty can substantially erode, or even eliminate, the net removal benefit. AIDC waste heat directly addresses this bottleneck by reducing the additional energy required for absorbent regeneration. The most important consequence is the observed feasibility flip: regions where standalone DAC struggles to achieve net removal can become carbon-negative when regeneration heat is supplied by AIDC waste heat. This result is particularly relevant for near-term deployment, because the power grid will not be fully decarbonized everywhere before the next wave of AI infrastructure expansion occurs.

The removal ratio provides a direct measure of whether DAC--AIDC integration can compensate for the operational emissions of AI computing. A ratio above unity is especially important because it indicates that the integrated system removes more CO$_2$ than the associated AIDC emits through electricity consumption. In 2024, several regions already approach or exceed this threshold. By 2030, under a future with more GPU-intensive AIDCs and cleaner electricity systems, several states achieve removal ratios far above 1. These values suggest that DAC--AIDC systems in favorable regions could offset not only their own operational emissions but also emissions from additional AI infrastructure. This shifts the role of integration from facility-level mitigation to regional carbon-removal planning.

These findings have implications for the siting of future AI data centers. Current siting decisions are often dominated by electricity price, land availability, water access, tax incentives, fiber connectivity, and grid interconnection capacity. Our results suggest that carbon-removal potential should also be considered. Regions with clean or rapidly decarbonizing grids, favorable climate conditions for DAC, and large recoverable waste-heat streams are better positioned to host AI infrastructure with lower net climate impact. In such regions, additional AIDC deployment may create stronger opportunities for coupling computing growth with carbon removal. Conversely, placing large AIDC clusters in regions with carbon-intensive grids and unfavorable DAC operating conditions may lock in higher emissions and reduce the effectiveness of future mitigation measures.

Overall, this study reframes the sustainability challenge of AI infrastructure. The growth of AIDCs will increase electricity demand and associated emissions, but it will also create concentrated, continuous, and increasingly high-quality waste-heat streams. If strategically integrated with DAC and deployed in regions with favorable grid and climate conditions, AIDCs could evolve from carbon-intensive digital infrastructure into anchors for distributed carbon-removal systems. The decarbonization of AI should therefore not rely only on procuring cleaner electricity; it should also exploit the thermodynamic characteristics of AI computing itself.


\section*{Methods}

\subsection*{Overview of the integrated system model}

We develop an integrated DAC--AIDC model in MATLAB, with the heat-pump cycle simulated in Aspen HYSYS. The model couples three modules: an AIDC cooling system, a vapor-compression heat pump, and a Lewatit VP OC 1065 / TVSA direct air capture (DAC) system. The heat pump consists of an evaporator, condenser, expansion valve, and compressor, using R245fa as the working fluid. Its coefficient of performance $\mathrm{COP}$ and heat-exchanger efficiency $\eta$ are parameterized as functions of the evaporation and condensation temperatures. The integrated model apply one month as time step and more details are provided in the Supplementary Information.

The model is indexed over months $t \in \mathcal{T} = \{1,\ldots,12\}$, U.S. states $s \in \mathcal{S}$, counties $c \in \mathcal{C}_{s}$ within each parent state, AI server classes $y \in \mathcal{Y} = \{\text{AI-8GPU},\,\text{AI-4GPU},\,\text{AI-2GPU},\,\text{AI Non-accelerated}\}$, and server load conditions $x \in \mathcal{X} = \{\text{train},\,\text{inference},\,\text{idle}\}$.

\subsection*{Waste heat allocation}

Following the 2024 Lawrence Berkeley National Laboratory report~\cite{shehabi20242024}, the installed AI server stock is resolved by server class $y$, rated power $P^{y}_{x,t}$ (W), and operating duration $\tau^{y}_{x,t}$ (h) under each load condition $x$ and month $t$. The annual national AI server electricity consumption is
\begin{equation}
E^{\mathrm{ser}}
\;=\;
\sum_{t \in \mathcal{T}}\, \sum_{y \in \mathcal{Y}}\, \sum_{x \in \mathcal{X}}
P^{y}_{x,t}\,\tau^{y}_{x,t},
\end{equation}
where $E^{\mathrm{ser}}$ (MWh) denotes the annual national AI server electricity consumption. Waste-heat recovery factors $\rho_{y}$ are used to represent the recoverable heat fraction for different server classes. These factors are determined by the share of heat dissipated through liquid-cooled components, including CPUs/GPUs and memory chips operating above $60\,^\circ\mathrm{C}$, that can be upgraded by the heat pump. The values of $\rho_{y}$ are 0.89, 0.78, 0.69, and 0.46 for AI-8GPU, AI-4GPU, AI-2GPU, and AI Non-accelerated servers, respectively~\cite{ahmed2021review,nvidia2025dgxsuperpod}. The total recoverable AIDC waste heat in month $t$ is
\begin{equation}
Q^{\mathrm{tot}}_{t}
\;=\;
\sum_{y \in \mathcal{Y}} \rho_{y}
\sum_{x \in \mathcal{X}} P^{y}_{x,t}\,\tau^{y}_{x,t},
\end{equation}
where $Q^{\mathrm{tot}}_{t}$ (MWh) is the national recoverable AIDC waste heat in month $t$. The recoverable heat is then spatially allocated using EPRI's 2023 state-level AIDC electricity profile~\cite{intelligence2024analyzing}. AI servers are first distributed across the top thirty DC states according to each state share $\alpha_{s}$ of national DC electricity use, with $\sum_{s \in \mathcal{S}} \alpha_{s} = 1$. The low-grade waste heat available in state $s$ and month $t$ is
\begin{equation}
Q^{\mathrm{lo}}_{s,t} \;=\; \alpha_{s}\, Q^{\mathrm{tot}}_{t}.
\end{equation}

State-level waste heat is further allocated to counties using DOE ``Speed to Power'' county-level DC capacity data~\cite{DOE_SpeedToPower_2025}. This is calculated by each county share in corresponding state $\gamma_{c}$, with $\sum_{c \in \mathcal{C}_{s}} \gamma_{c} = 1$. The low-grade waste heat available in county $c \in \mathcal{C}_{s}$ and month $t$ is
\begin{equation}
Q^{\mathrm{lo}}_{s,c,t} \;=\; \gamma_{c}\, Q^{\mathrm{lo}}_{s,t}.
\end{equation}

\subsection*{Thermal coupling and DAC capture}

The heat pump upgrades $Q^{\mathrm{lo}}_{s,c,t}$ from the AIDC warm return temperature ($\sim\!60\,^\circ\mathrm{C}$) to the DAC regeneration set point ($100\,^\circ\mathrm{C}$). Heat transfer is considered across both the evaporator and condenser, with constant heat exchanger efficiency $\eta$. COP is the coefficient of performance of heat pump, which is different in the DAC-AIDC system using high-temperature heat pump and the standalone DAC system using supercritical heat pump. The standard heat pump relation $Q^{\mathrm{hi}} = Q^{\mathrm{lo}}\cdot \mathrm{COP}/(\mathrm{COP}-1)$ i
s also applied. The high-grade heat delivered for DAC desorption is therefore
\begin{equation}
Q^{\mathrm{hi}}_{s,c,t}
\;=\;
Q^{\mathrm{lo}}_{s,c,t}\,\eta^{2}\,
\frac{\mathrm{COP}}{\mathrm{COP} - 1},
\end{equation}
where $Q^{\mathrm{hi}}_{s,c,t}$ (MWh) is the regeneration heat available in county $c$ and month $t$. $q^{\mathrm{th}}_{s,c,t}$ (MWh/t) and $E_{s,c,t}$ (MWh/t) denote the per-tonne thermal and electrical demand of the DAC operational cycle, respectively. Both values are evaluated from county-level monthly temperature and relative humidity, as described in the SI. Modular DAC capacity is assumed to be sized to fully use the available high-grade heat. The total monthly CO$_2$ capture in county $c$ is then
\begin{equation}
m^{\mathrm{tot}}_{s,c,t}
\;=\;
\frac{Q^{\mathrm{hi}}_{s,c,t}}{q^{\mathrm{th}}_{s,c,t}},
\qquad
m^{\mathrm{tot}}_{s,c}
\;=\;
\sum_{t \in \mathcal{T}} m^{\mathrm{tot}}_{s,c,t},
\end{equation}
where $m^{\mathrm{tot}}_{s,c,t}$ (t) is the gross monthly capture and $m^{\mathrm{tot}}_{s,c}$ is the corresponding annual gross capture in county $c$. Additionally, the total grid electricity demand consists of heat pump $E^{\mathrm{hp}}_{s,c,t}$ and DAC fan/auxiliary electricity use $E^{\mathrm{dac}}_{s,c,t}$, both in MWh:
\begin{equation}
E^{\mathrm{tot}}_{s,c,t}
\;=\;
E^{\mathrm{hp}}_{s,c,t} + E^{\mathrm{dac}}_{s,c,t}.
\end{equation}
Net annual CO$_2$ removal is calculated by subtracting the indirect emissions from grid electricity use. These emissions are evaluated using the state grid carbon emission factor (CEF) $\mu_{s}$ (kg\,CO$_2$/MWh):
\begin{equation}
m^{\mathrm{net}}_{s,c}
\;=\;
\sum_{t \in \mathcal{T}}
\bigl( m^{\mathrm{tot}}_{s,c,t} \;-\; \mu_{s}\, E^{\mathrm{tot}}_{s,c,t} \bigr).
\end{equation}

\subsection*{Removal ratio}

The carbon-removal benefit of DAC--AIDC thermal integration is evaluated from the AIDC perspective. The \emph{removal ratio} is defined as the ratio between net CO$_2$ removal by DAC and the indirect CO$_2$ emissions from AIDC electricity consumption:
\begin{equation}
\varphi_{s,c}
\;=\;
\frac{m^{\mathrm{net}}_{s,c}}
     {\mu_{s}\,\alpha_{s}\,\gamma_{c}\,\lambda\, E^{\mathrm{ser}}}.
\end{equation}
The denominator $\mu_{s}\,\alpha_{s}\,\gamma_{c}\,\lambda\,E^{\mathrm{ser}}$ represents the annual indirect CO$_2$ emissions of the AIDC located in county $c$ due to grid electricity use. $\mathrm{\lambda}$ refers to the power usage effectiveness (PUE), which is applied to calculate the total electricity consumption of the AIDC. A value of $\varphi_{s,c} = 1$ indicates net-zero AIDC operation. A value of $\varphi_{s,c} > 1$ indicates that the integrated system offsets all indirect AIDC emissions and removes additional CO$_2$. When $\varphi_{s,c} < 0$, the DAC system is not suitable to deploy as its operation causes additional CO$_2$ emissions.
\subsection*{Economic model}

The levelized cost of CO$_2$ capture (LCCC) is calculated at the state level by combining per-tonne operating expenditure with annualized capital cost. For a state electricity price $\pi_{s}$ (\$/MWh) and monthly weighting factor $\beta_{t}$ ($\sum_{t \in \mathcal{T}} \beta_{t} = 1$, refers to annual COCO$_2$ capture occurring in month t), the operating expenditure is
\begin{equation}
\mathrm{OPEX}_{s}
\;=\;
\pi_{s}
\sum_{t \in \mathcal{T}}
\beta_{t}
\!\left(
E_{s,t}
\;+\;
\frac{q^{\mathrm{th}}_{s,t}\,\eta}{\mathrm{COP}}
\right),
\end{equation}
where $\mathrm{OPEX}_{s}$ (\$/t) is the state-level operating expenditure per tonne CO$_2$ captured. The LCCC is then
\begin{equation}
\mathrm{LCCC}_{s}
\;=\;
\mathrm{OPEX}_{s} \;+\; \sum_{e \in \mathcal{E}} \mathrm{CAPEX}_{s,e},
\end{equation}
where $\mathrm{CAPEX}_{s,e}$ (\$/t) is the annualized per-tonne capital cost of equipment item $e \in \mathcal{E}$, including heat-pump components, DAC modules, fans, and absorbent, in state $s$. More details are provided in SI.




\bibliographystyle{unsrturl}
\bibliography{References.bib}

@article{ebrahimi2014review,
  title={A review of data center cooling technology, operating conditions and the corresponding low-grade waste heat recovery opportunities},
  author={Ebrahimi, Khosrow and Jones, Gerard F and Fleischer, Amy S},
  journal={Renewable and sustainable energy reviews},
  volume={31},
  pages={622--638},
  year={2014},
  publisher={Elsevier}
}

@article{shehabi20242024,
  title={2024 united states data center energy usage report},
  author={Shehabi, Arman and Newkirk, Alex and Smith, Sarah J and Hubbard, Alex and Lei, Nuoa and Siddik, Md Abu Bakar and Holecek, Billie and Koomey, Jonathan and Masanet, Eric and Sartor, Dale},
  year={2024}
}

@article{intelligence2024analyzing,
  title={Analyzing Artificial Intelligence and Data Center Energy Consumption},
  author={Intelligence, Powering},
  journal={Electric Power Research Institute: Washington, DC, USA},
  year={2024}
}

@techreport{IEA2025_energy_ai,
  author       = {{International Energy Agency}},
  title        = {Energy and AI},
  institution  = {IEA},
  address      = {Paris},
  year         = {2025},
  url          = {https://www.iea.org/reports/energy-and-ai}
}

@article{azarifar2024liquid,
  title={Liquid cooling of data centers: A necessity facing challenges},
  author={Azarifar, Mohammad and Arik, Mehmet and Chang, Je-Young},
  journal={Applied Thermal Engineering},
  volume={247},
  pages={123112},
  year={2024},
  publisher={Elsevier}
}

@online{mckinsey2024_ai_power,
  author       = {{McKinsey \& Company}},
  title        = {AI Power: Expanding Data Center Capacity to Meet Growing Demand},
  year         = {2024},
  month        = oct,
  url          = {https://www.mckinsey.com/industries/technology-media-and-telecommunications/our-insights/ai-power-expanding-data-center-capacity-to-meet-growing-demand}
}

@article{guidi2024environmental,
  title={Environmental burden of United States data centers in the artificial intelligence era},
  author={Guidi, Gianluca and Dominici, Francesca and Gilmour, Jonathan and Butler, Kevin and Bell, Eric and Delaney, Scott and Bargagli-Stoffi, Falco J},
  journal={arXiv preprint arXiv:2411.09786},
  year={2024}
}

@article{xiao2025environmental,
  title={Environmental impact and net-zero pathways for sustainable artificial intelligence servers in the USA},
  author={Xiao, Tianqi and Nerini, Francesco Fuso and Matthews, H Damon and Tavoni, Massimo and You, Fengqi},
  journal={Nature Sustainability},
  pages={1--13},
  year={2025},
  publisher={Nature Publishing Group UK London}
}

@techreport{epa2023_passenger_vehicle_emissions,
  author       = {{U.S. Environmental Protection Agency}},
  title        = {Greenhouse Gas Emissions from a Typical Passenger Vehicle},
  institution  = {Office of Transportation and Air Quality, U.S. Environmental Protection Agency},
  number       = {EPA-420-F-23-014},
  year         = {2023},
  month        = jun,
  url          = {https://www.epa.gov/greenvehicles/greenhouse-gas-emissions-typical-passenger-vehicle}
}

@article{patterson2021carbon,
  title={Carbon emissions and large neural network training},
  author={Patterson, David and Gonzalez, Joseph and Le, Quoc and Liang, Chen and Munguia, Lluis-Miquel and Rothchild, Daniel and So, David and Texier, Maud and Dean, Jeff},
  journal={arXiv preprint arXiv:2104.10350},
  year={2021}
}

@inproceedings{zhang2023global,
  title={The global energy impact of raising the space temperature for high-temperature data centers, Cell Rep},
  author={Zhang, Y and Li, H and Wang, S},
  booktitle={Phys. Sci},
  volume={4},
  number={10},
  year={2023}
}

@article{tervo2025reducing,
  title={Reducing district heating carbon dioxide emissions with data center waste heat--Region perspective},
  author={Tervo, Seela and Syri, Sanna and Hiltunen, Pauli},
  journal={Renewable and Sustainable Energy Reviews},
  volume={208},
  pages={114992},
  year={2025},
  publisher={Elsevier}
}

@article{townsend2024our,
  title={Our first offsite heat recovery project lands in Finland},
  author={Townsend, Ben},
  journal={Google},
  year={2024}
}

@article{santin2020feasibility,
  title={Feasibility limits of using low-grade industrial waste heat in symbiotic district heating and cooling networks},
  author={Santin, Maurizio and Chinese, Damiana and De Angelis, Alessandra and Biberacher, Markus},
  journal={Clean Technologies and Environmental Policy},
  volume={22},
  number={6},
  pages={1339--1357},
  year={2020},
  publisher={Springer}
}

@article{jing2026large,
  title={Large-scale Long-distance Data Center Waste Heat District Heating System: Exergy Analysis, System Performance, and Application},
  author={Jing, Yang and Zhan, Ziyang and Li, Min and Cai, Rongbiao and Xie, Xiaoyun and Jiang, Yi},
  journal={Energy},
  pages={140076},
  year={2026},
  publisher={Elsevier}
}

@article{keith2018process,
  title={A process for capturing CO2 from the atmosphere},
  author={Keith, David W and Holmes, Geoffrey and Angelo, David St and Heidel, Kenton},
  journal={Joule},
  volume={2},
  number={8},
  pages={1573--1594},
  year={2018},
  publisher={Elsevier}
}

@article{mcqueen2021review,
  title={A review of direct air capture (DAC): scaling up commercial technologies and innovating for the future},
  author={McQueen, Noah and Gomes, Katherine Vaz and McCormick, Colin and Blumanthal, Katherine and Pisciotta, Maxwell and Wilcox, Jennifer},
  journal={Progress in Energy},
  volume={3},
  number={3},
  pages={032001},
  year={2021},
  publisher={IOP Publishing}
}

@book{international2022direct,
  title={Direct Air Capture: A Key Technology for Net Zero},
  author={International Energy Agency},
  year={2022},
  publisher={OECD Publishing}
}

@techreport{zhou2021low_temp_dac,
  author       = {Zhou, Xiaoliang and others},
  title        = {Low Regeneration Temperature Sorbent for Direct Air Capture of CO$_2$},
  institution  = {U.S. Department of Energy, National Energy Technology Laboratory (NETL)},
  year         = {2021},
  url          = {https://netl.doe.gov/sites/default/files/netl-file/21CMOG_CDRR_Zhou.pdf}
}

@article{sabatino2021comparative,
  title={A comparative energy and costs assessment and optimization for direct air capture technologies},
  author={Sabatino, Francesco and Grimm, Alexa and Gallucci, Fausto and van Sint Annaland, Martin and Kramer, Gert Jan and Gazzani, Matteo},
  journal={Joule},
  volume={5},
  number={8},
  pages={2047--2076},
  year={2021},
  publisher={Elsevier}
}

@article{sievert2024considering,
  title={Considering technology characteristics to project future costs of direct air capture},
  author={Sievert, Katrin and Schmidt, Tobias S and Steffen, Bjarne},
  journal={Joule},
  volume={8},
  number={4},
  pages={979--999},
  year={2024},
  publisher={Elsevier}
}

@article{wiegner2022optimal,
  title={Optimal design and operation of solid sorbent direct air capture processes at varying ambient conditions},
  author={Wiegner, Jan F and Grimm, Alexa and Weimann, Lukas and Gazzani, Matteo},
  journal={Industrial \& Engineering Chemistry Research},
  volume={61},
  number={34},
  pages={12649--12667},
  year={2022},
  publisher={ACS Publications}
}

@article{leonzio2022innovative,
  title={Innovative process integrating air source heat pumps and direct air capture processes},
  author={Leonzio, Grazia and Shah, Nilay},
  journal={Industrial \& Engineering Chemistry Research},
  volume={61},
  number={35},
  pages={13221--13230},
  year={2022},
  publisher={ACS Publications}
}

@article{ge2024innovative,
  title={Innovative process integrating high temperature heat pump and direct air capture},
  author={Ge, Bingyao and Zhang, Man and Hu, Bin and Wu, Di and Zhu, Xuancan and Eicker, Ursula and Wang, Ruzhu},
  journal={Applied Energy},
  volume={355},
  pages={122229},
  year={2024},
  publisher={Elsevier}
}

@article{odeh2025techno,
  title={Techno-economic assessment of waste heat-powered direct air capture in the refinery and petrochemical sectors in Saudi Arabia},
  author={Odeh, Naser and Apeaning, Raphael W and Rowaihy, Feras},
  journal={Carbon Capture Science \& Technology},
  volume={16},
  pages={100451},
  year={2025},
  publisher={Elsevier}
}

@misc{DOE_SpeedToPower_2025,
  author       = {{U.S. Department of Energy}},
  title        = {Speed to Power Initiative},
  year         = {2025},
  howpublished = {\url{https://www.energy.gov/speed-to-power}}
}

@article{diaz2025flipping,
  title={Flipping the switch: carbon-negative and water-positive data centers through waste heat utilization},
  author={D{\'\i}az-Mar{\'\i}n, Carlos D and Berquist, Zachary J},
  journal={Energy \& Environmental Science},
  volume={18},
  number={18},
  pages={8403--8413},
  year={2025},
  publisher={Royal Society of Chemistry}
}

@article{masanet2020recalibrating,
  title={Recalibrating global data center energy-use estimates},
  author={Masanet, Eric and Shehabi, Arman and Lei, Nuoa and Smith, Sarah and Koomey, Jonathan},
  journal={Science},
  volume={367},
  number={6481},
  pages={984--986},
  year={2020},
  publisher={American Association for the Advancement of Science}
}

@techreport{fout2022direct,
  title={Direct Air Capture Case Studies},
  author={Fout, Timothy and Zoelle, Alex and Homsy, Sally and Valentine, Jessica and Roy, Naksha and Kilstofte, Aaron and Sturdivan, Mike and Steutermann, Mark and Woods, Mark},
  year={2022},
  institution={National Energy Technology Laboratory (NETL), Pittsburgh, PA, Morgantown, WV~…}
}

@inproceedings{patel2024techno,
  title={Techno-economic analysis of sorbent-based direct air capture informed by EPC input and recent technological advancements},
  author={Patel, Kshitij and Teel, Troy and Henry, Samuel and Mantripragada, Hari and Zoelle, Alexander and Fout, Timothy and Homsy, Sally},
  booktitle={Proceedings of the 17th Greenhouse Gas Control Technologies Conference (GHGT-17)},
  pages={20--24},
  year={2024}
}

@article{sendi2022geospatial,
  title={Geospatial analysis of regional climate impacts to accelerate cost-efficient direct air capture deployment},
  author={Sendi, Marwan and Bui, Mai and Mac Dowell, Niall and Fennell, Paul},
  journal={One Earth},
  volume={5},
  number={10},
  pages={1153--1164},
  year={2022},
  publisher={Elsevier}
}

@article{d2024integrating,
  title={Integrating direct air capture with algal biofuel production to reduce cost, energy, and GHG emissions},
  author={D’Souza, Shavonn and Johnston, Jaden and Thomas, VM and Harris, Kylee and Tan, Eric CD and Chance, Ronald R and Yuan, Yanyui},
  journal={Journal of CO2 Utilization},
  volume={86},
  pages={102911},
  year={2024},
  publisher={Elsevier}
}

@article{de2023growing,
  title={The growing energy footprint of artificial intelligence},
  author={De Vries, Alex},
  journal={Joule},
  volume={7},
  number={10},
  pages={2191--2194},
  year={2023},
  publisher={Elsevier}
}

@misc{nvidia2025dgxsuperpod,
  title        = {NVIDIA DGX SuperPOD: Data Center Design Featuring NVIDIA DGX H100 Systems},
  author       = {{NVIDIA Corporation}},
  year         = {2025},
  note         = {Planning a Data Center Deployment, Power and Heat Dissipation section.},
  howpublished = {https://docs.nvidia.com/dgx-superpod/design-guides/dgx-superpod-data-center-design-h100/latest/planning.html}
}

@article{ahmed2021review,
  title={A review of data centers energy consumption and reliability modeling},
  author={Ahmed, Kazi Main Uddin and Bollen, Math HJ and Alvarez, Manuel},
  journal={IEEE access},
  volume={9},
  pages={152536--152563},
  year={2021},
  publisher={IEEE}
}

@techreport{gagnon2025cambium,
  title       = {Cambium 2024 Scenario Descriptions and Documentation},
  author      = {Gagnon, Pieter and Sanchez Perez, Pedro Andres and Florez, Julian and Morris, James and Llerena Velasquez, Marck and Eisenman, Jordan},
  institution = {National Renewable Energy Laboratory},
  address     = {Golden, CO},
  number      = {NREL/TP-6A40-93005},
  year        = {2025},
  month       = apr,
}

@article{alissa2025using,
  title={Using life cycle assessment to drive innovation for sustainable cool clouds},
  author={Alissa, Husam and Nick, Teresa and Raniwala, Ashish and Arribas Herranz, Alberto and Frost, Kali and Manousakis, Ioannis and Lio, Kari and Warrier, Brijesh and Oruganti, Vaidehi and DiCaprio, TJ and others},
  journal={Nature},
  volume={641},
  number={8062},
  pages={331--338},
  year={2025},
  publisher={Nature Publishing Group UK London}
}

@misc{DellOroGroup2026LiquidCoolingMarket,
  author       = {{Dell'Oro Group}},
  title        = {Data Center Liquid Cooling Market to Approach \$7 Billion by 2029 as {AI} Deployments Accelerate, According to {Dell'Oro Group}},
  year         = {2026},
  month        = jan,
  howpublished = {\url{https://www.delloro.com/news/data-center-liquid-cooling-market-to-approach-7-billion-by-2029-as-ai-deployments-accelerate/}},
  note         = {Press release, accessed 2026-05-12}
}

@misc{Davis2024WaterColdPlatesDLC,
  author       = {Davis, Jacqueline},
  title        = {Water Cold Plates Lead in the Small, but Growing, World of {DLC}},
  year         = {2024},
  month        = oct,
  howpublished = {\url{https://journal.uptimeinstitute.com/water-cold-plates-lead-in-the-small-but-growing-world-of-dlc/}},
  note         = {Uptime Institute Journal, accessed 2026-05-12}
}

\section*{Acknowledgments}

\section*{Author Contributions}

\section*{Declaration of Competing Interest} 

\end{document}